\documentclass[10pt]{article}

\usepackage{amssymb,amsmath}
\usepackage{amsthm}
\usepackage{graphicx}
\usepackage{pdfsync}
\usepackage{epstopdf}

\usepackage{color}

 \newtheorem{thm}{Theorem}[section]
 \newtheorem{defn}[thm]{Definition}
 \newtheorem{ex}[thm]{Example}
 \newtheorem{cor}[thm]{Corollary}
 \newtheorem{lem}[thm]{Lemma}
 \newtheorem{prop}[thm]{Proposition}

\newtheorem{obs}[thm]{Observation}

\newtheorem{alg}[thm]{Algorithm}

\def\G{{\mathcal G}}

\def\S{{\mathcal S}}

\def\naturals{{\mathbb N}}

\def\reals{{\mathbb R}}
\def\complex{{\mathbb C}}

\def\rmnmats #1 #2{M_{#1 #2}\,({\mathbb R})}
\def\cmnmats #1 #2{M_{#1 #2}\,({\mathbb C})}

\def\mnmats #1 #2{M_{#1 #2}}

\def\matdim #1 #2{#1 \times #2}

\def\rk{{\rm rank}}

\newcommand{\sS}{\mathcal{S}^-}
\newcommand{\smr}{\operatorname{mr}^-}
\newcommand{\srs}{\operatorname{r}^-}
\newcommand{\sMR}{\operatorname{MR}^-} 

\newcommand{\mr}{\operatorname{mr}}

\newcommand{\sZ}{\operatorname{Z}^-}

\newcommand{\match}{\operatorname{match}}

\newcommand{\dm}{\operatorname{diam}}

\newcommand{\x}{\times}
\newcommand{\bit}{\begin{itemize}}
\newcommand{\eit}{\end{itemize}}
\newcommand{\ben}{\begin{enumerate}}
\newcommand{\een}{\end{enumerate}}
\newcommand{\beq}{\begin{equation}}
\newcommand{\eeq}{\end{equation}}
\newcommand{\bea}{\begin{eqnarray*}}
\newcommand{\eea}{\end{eqnarray*}}
\newcommand{\bpf}{\begin{proof}}
\newcommand{\epf}{\end{proof}}

\title{Acyclic and unicyclic graphs whose minimum skew rank is equal to the minimum skew rank of a diametrical path}

\author{Luz~M.~DeAlba\thanks{Department of Mathematics and Computer
Science, Drake University, Des Moines, IA 50311, USA
(luz.dealba@drake.edu).}}

\begin{document}

\maketitle

\bigskip

\begin{abstract}
The minimum skew rank of a simple graph $G$ over the field of real numbers, $\smr (G)$, is  the
smallest possible rank among all real skew-symmetric  matrices  whose $(i,j)$-entry (for $i\neq j$) is nonzero whenever   $\{i,j\}$ is an edge in $G$ and is zero otherwise.  In this paper we give an algorithm for computing the minimum skew rank of a connected unicyclic graph, and classify all connected acyclic and connected unicyclic graphs $G$, for which $\smr (G) = \smr (P)$, where $P$ is a diametrical path of $G$.
\end{abstract}



\section{Introduction}

A {\em simple undirected graph} is a pair  $G = (V_G, E_G)$, where $V_G$ is the (finite, nonempty) set of vertices of $G$ and $E_G$ is the set of edges, where an edge  is a two-element subset of vertices. In this paper, all graphs are simple, that is, no loops or multiple edges are allowed. The {\em order of a graph} $G$, denoted $|G|$, is the number of vertices of $G$. A matrix $A \in F^{n \x n}$ ($F$ a field) is {\em skew-symmetric} (respectively, {\em symmetric}) if $A^T = -A$ (respectively, $A^T = A$); for $A \in \complex^{n \x n}$, $A$ is {\em Hermitian} if $A^*= A$, where $A^*$ denotes the conjugate transpose of $A$. Unless explicitly stated otherwise, all symmetric and skew-symmetric matrices are assumed to be real.

For an $n \x n$ symmetric, skew-symmetric or Hermitian matrix $A$, the {\em graph of $A$}, denoted $\G(A)$, is the graph with vertices $\{v_1, . . . , v_n \}$ and edges $\{ \{v_i, v_j \} : a_{ij} \ne 0, 1 \le i < j \le n \}$. Note that for symmetric and Hermitian matrices the diagonal is ignored in determining $\G(A)$, but for a skew-symmetric matrix all diagonal entries must be 0. The classic minimum rank problem involves symmetric matrices and has been studied extensively. See~\cite{07FH} for a survey of known results and discussion of the motivation for the minimum rank problem; an extensive bibliography is also provided there. The study of minimum skew rank of matrices began in the recent article~\cite{10IMA}, and is an area where many open problems exist.

In this paper we give an algorithm for computing the minimum skew rank of a connected unicyclic graph, and classify all connected acyclic and connected unicyclic graphs $G$, for which $\smr (G) = \smr (P)$, where $P$ is a diametrical path of $G$. These problems were first studied for symmetric matrices in~\cite{07FH} and in~\cite{08BFS}. In Section~\ref{notation} we introduce necessary notation and terminology. In Section~\ref{prelim} we provide preliminary minimum skew rank results. In Section~\ref{tree} we give our main results on trees, and in Section~\ref{cycle} we give our main results on connected unicyclic graphs. Section~\ref{examples} contains examples of connected graphs with many cycles for which the minimum skew rank is equal to the minimum skew rank of a diametrical path.

\section{Notation and Terminology}\label{notation}

For a graph $G$, we let $\S(G) = \{ A \in \reals^{n \x n} : A^T = A, \G(A) = G \}$ and $\sS(G) = \{ A \in \reals^{n \x n} : A^T = -A, \G(A) = G \}$ denote, respectively,  the set of real symmetric matrices, and  the set of real skew-symmetric matrices described by  $G$. We let $\mr(G) = \min \{ \rk(A) : A \in \S(G) \}$, $\smr(G) = \min \{ \rk(A) : A \in \sS(G) \}$, and $\sMR(G) = \max \{ \rk(A) : A \in \sS(G) \}$ denote, respectively, the {\em minimum rank}, the {\em minimum skew rank}, and the {\em maximum skew rank} of $G$.

For a graph $G$, the {\em degree of a vertex $v \in V_G$}, denoted by $\deg_G (v)$, is the number of edges adjacent to $v$. A {\em path} on $n$ vertices is denoted by $P_n$, and its vertices numbered sequentially $v_1, v_n, \dots, v_{n-1}, v_n$. The {\em length of a path $P$} is the number of edges in the path, which is equal to  $|P|-1$. The {\em distance}  between two vertices $u$ and $v$, in a connected graph $G$, is the length of the shortest path joining $u$ and $v$, and is denoted by $d_G (u, v)$. The {\em diameter} of a connected graph $G$ is given by    $\dm (G) = \max \{ d_G (u, v) :  u, v \in V_G \}$. An induced path in a graph $G$, of length $\dm (G)$ is referred to as a {\em diametrical path} of $G$. If $P$ is an induced path of the longest possible length in a graph $G$, then $\dm (G) \le |P|-1$. For a tree $T$  equality holds, that is, $\dm (T)= |P|-1$.

A tree $T$ is called a {\em centipede} if it consists of a path $P$ with additional edges, called {\em legs} attached to the nonpendant vertices, these vertices with legs are {\em joints} of $P$.  A centipede is {\em regular} if no two consecutive vertices in $P$
are joints, and irregular otherwise.

A {\em matching} in a graph $G$ is a set of edges $M = \{ \{v_{ i_1}, v_{j_1}\}, \{ v_{i_2}, v_{j_2} \}, \dots, \{ v_{i_k}, v_{j_k} \} \}$ such that no endpoints are shared. The vertices in $V_G$ that are adjacent to the edges in $M$ are called {\em $M$-saturated} vertices, and those vertices in $V_G$ that are not adjacent to the edges in $M$ are called {\em $M$-unsaturated} vertices. A {\em perfect matching} in a graph $G$ is a matching that saturates all vertices of $G$. A {\em maximum matching} in a graph $G$ is a matching of maximum order among all matchings in $G$. The {\em matching number} of a graph $G$, denoted by $\match(G)$ is the number of edges in a maximum matching. A path in a graph $G$ is called {\em $M$-alternating} if it alternates between edges in $M$ and edges not in $M$. An $M$-alternating path is called {\em $M$-augmenting} whenever its endpoints are $M$-unsaturated.

\section{Preliminary Minimum Skew Rank Results}\label{prelim}

This section contains results that we use in Sections~\ref{tree} and~\ref{cycle} in our main results, and extend some known results on minimum rank to skew minimum rank.

Part~\ref{zf}, in Observation~\ref{mrlesssmr}, follows from Proposition~4.2 in \cite{AIM}, and Propositions~3.5 and~3.7 in~\cite{10IMA}. The zero forcing number, $ Z (T)$, is defined in ~\cite{AIM},  and the skew zero forcing number, $\sZ (T)$, is defined in~\cite{10IMA}. Part~\ref{dmlemsr} follows from Part~\ref{zf}.

\begin{obs}\label{mrlesssmr}
\ben
\item\label{induced}\cite[Observation~1.6]{10IMA} If $H$ is an induced subgraph of $G$, then $\smr (H) \le \smr (G)$.

\item\label{zf} If $T$ is a tree, then $\mr (T) = |T| - Z (T) \le |T| - \sZ (T) \le \smr (T)$.
\item\label{dmlemsr} For a graph $G$ with diametrical path $P$, $\dm (G) \le \smr(P)  \le \smr (G)$.
\een
\end{obs}

\begin{thm}\label{berge} (\cite[p. 109]{01W})
A matching $M$ in a graph $G$ is a maximum matching in $G$ if and only if $G$ has no $M$-augmenting path.
\end{thm}

\begin{thm}\label{MR(G)}~\cite[Theorem~2.5]{10IMA} For a graph $G$, $\sMR (G) = 2 \match(G)$, and every even rank between $\smr (G)$ and $\sMR (G)$ is realized by a matrix in $\sS (G)$.
\end{thm}

\begin{thm}\label{upf}~\cite[Theorem~2.6]{10IMA} For a graph $G$, $\smr(G) = |G| =
\sMR(G)$ if and only if $G$ has a unique perfect matching.
\end{thm}

\begin{thm}\label{smrtree}\cite[Theorem~2.7]{10IMA}
If $T$ is a tree, then $\smr(T) = 2 \match(T) = \sMR(T)$.
\end{thm}

\begin{prop}\label{smr(P)}~\cite[Proposition 4.1]{10IMA}
For a path $P_n$ on $n$ vertices,
$$\smr (P_n)
 = \left\{
\begin{array}{lll}
n & \mbox{if}  & n \mbox{ is even}\\
n - 1 & \mbox{if} & n \mbox{ is odd}
\end{array}
\right..$$
\end{prop}

\begin{prop}\label{smr(C)}~\cite[Proposition 4.2]{10IMA}
For a cycle $C_n$ on $n$ vertices,
$$\smr (C_n)
 = \left\{
\begin{array}{lll}
n - 2 & \mbox{if}  &n \mbox{ is even}\\
n - 1 & \mbox{if} & n \mbox{ is odd}
\end{array}
\right..$$
\end{prop}

Lemma~\ref{fallat} extends the results of Lemma~4 of~\cite{08BFS} to minimum skew rank;  the condition $\smr (G) = \dm (G)$ in~~\cite{08BFS} is replaced with $\smr (G) = \smr (P)$, where $P$ is a diametrical path of $G$.

\begin{lem}\label{fallat} Let $P$ be a diametrical path of a connected graph $G$, where $\smr (G) = \smr (P)$, and let $H$ be a set of vertices that do not lie in $P$. Then $\smr (G-H) = \smr (G)$.
\end{lem}

\section{Trees}\label{tree}

For $n$ odd, $P_n$ has exactly $\frac{n+1}{2}$ distinct maximum matchings, each leaving exactly one of the odd-numbered vertices unsaturated. When $n$ is even, $P_n$ has a unique perfect matching.

\begin{lem}\label{centi} Let $T$ be a tree with diametrical path $P$.  If $\smr (T) = \smr (P)$, then $T$ is a centipede. Moreover, if $|P|$ is odd, then $T$ is a regular centipede.
\end{lem}

\bpf
If $|P|$ is odd, by Observation~\ref{mrlesssmr}, $\dm (T) \le \mr (T) \le \smr (T) = \smr (P) = \dm (T)$, so  $T$ is a regular centipede([Theorem~6]~\cite{08BFS}). If $|P|$ is even and $T$ is not a centipede, the unique perfect matching in $P$ can be extended to a matching in $T$. By Theorem~\ref{smrtree}, $\smr (T) > \smr (P)$.
\epf
 
\begin{thm}\label{maintree} Let $T$ be a tree and $P_n$ a diametrical path of $T$.  If $n$ is odd $\smr (T) = \smr (P_n)$ if and only if $T$ is a regular centipede with the property that all joints are even-indexed vertices of $P_n$, and if $n$ is even $\smr (T) = \smr (P_n)$ if and only if $T$ is a centipede with the property that no pair of  joints $v_r$ and $v_s$ of  $P_n$ satisfy $r$ odd, $s$ even, and $3 \le r < s \le n-2$.
\end{thm}

\bpf
If $\smr (T) = \smr (P_n)$, by Lemma~\ref{centi}, $T$ is centipede (regular if $n$ is odd). If $n$ is odd and $v_i$, is a joint with $i$ odd, a maximum matching in $P_n$ can clearly be extended to a matching in $T$. If $n$ is even, $P_n$ has  joints $v_r$ and $v_s$,  $r$ odd, $s$ even, $3 \le r < s \le n-2$, and $w_1, w_2 \in V_T - V_{P_n}$ are adjacent to $v_r$ and $v_s$,  respectively, then the path $w_1 v_r  v_{r+1} \dots v_{s-1} v_s, w_2$ is an $M$-augmenting path in $T$, where $M$ is the unique perfect matching in $P$, in both cases, by Theorem~\ref{smrtree}, we have a contradiction.

Suppose $T$ is a centipede. If $T$ has at most one joint ($v_i$, $i$ even if $n$ is odd), it is easy to see that $\smr (T) = \smr (P_n)$, thus we assume that $T$ has at least two joints, $v_r$, and $v_s$, and let $M$ be a maximum matching in $T$.

Suppose $n$ is odd, and $r$ and $s$ are even. Without loss of generality, we can assume $s = r +2k$, for some $k \in \naturals$. Let $w_1, w_2 \in V_T - V_{P_n}$ be (clearly $M$-unsaturated) adjacent to $v_r$ and $v_s$,  respectively. The distance between $v_r$ and $v_s$ is even, thus $T$ cannot have an $M$-alternating path starting at $w_1$ and ending at $w_2$. Also, $T$ does not have an $M$-alternating path starting at either $w_1$ or $w_2$ and ending at the only $M$-unsaturated vertex on $P_n$.  It is also clear that there is no $M$-alternating path starting at a vertex adjacent to $v_r$ (or $v_s$), and ending at another vertex adjacent to $v_r$ (or $v_s$).

Suppose $n$ is even. If both $s$ and $r$ are odd or both are even, then the distance between joints is even, and $T$ cannot have an $M$-augmenting path starting at a vertex adjacent to $v_s$ and ending at a vertex adjacent to $v_r$. If $v_s$ and $v_r$ are joints, with $s$ even, $r$ odd, and $2 \le s < r \le n-1$, then  let  $w_1, w_2 \in V_T - V_{P_n}$ be adjacent to $v_s$ and $v_r$, respectively. Now $\{ v_{s-1},  v_s \} \in M$, and $\{ w_1, v_s \}, \{ v_s,  v_{s+1} \} \notin M$, which implies $T$ cannot have an $M$-augmenting path starting at $w_1$ and ending at $w_2$.

Since in both cases $T$ has no $M$-augmenting path, it follows from Theorem~\ref{berge} that $M$ is also a maximum matching in $T$, and  from Theorem~\ref{smrtree} we have $\smr (T) = \smr (P_n)$.
\epf

\section{Unicyclic Graphs}\label{cycle}

A vertex in a graph $G$ is a {\em cut-vertex} if the number of components of $G-v$ is greater than the number of components of $G$. The {\em skew rank-spread of a graph $G$ at a vertex $v$} is defined as $\srs_v (G) = \smr (G) - \smr (G-v).$ Note that $\srs_v (G)$ has a value of either 0 or 2.

\begin{thm}\label{cutvertex}~(\cite{09D}) Let $v$ be a cut-vertex of $G$. For $i = 1, 2, \dots, h$, let $W_i \subseteq V(G)$ be the vertices of the $i$-th component of $G-v$, and let $G_i$ be the subgraph of $G$ induced by $\{ v \} \cup W_i$. Then
$$\srs_v(G) = \max_{i=1, 2, \dots, h} \srs_v( G_i),$$
and
$$\smr (G) = \left\{
\begin{array}{ll}
\sum_{i=1}^h \smr (G_i - v) & \mbox{if } \srs_v (G_i) = 0 \mbox{ for all } i = 1, 2, \dots, h,\\
\sum_{i=1}^h \smr (G_i - v) + 2 & \mbox{if } \srs_v (G_i) = 2 \mbox{ for some } i = 1, 2, \dots, h.
\end{array}
\right.$$
\end{thm}

\begin{defn}\label{pndan}
Let $C_n$ be an $n$-cycle and let $V' = \left\{ v_1, v_2, \dots, v_k \right\} \subseteq V_{C_n}, 0 < k \le n$. The unicyclic graph $U$, obtained from $C_n$ by appending $p_i$ leaves, $1 \le p_i$, to each vertex $v_i \in  V'$ is called a {\em partial $n$-dandelion}. If $p_i = 1$, and $k < n$, $U$ is a {\em partial $n$-sun}, and if $p_i = 1$, and $k=n$, then $U$ is an {\em $n$-sun}.
\end{defn}

\begin{ex}\label{pnsun}
If $U$ is a partial $n$-dandelion with $|V'|=k$, then $\smr (U) = 2 k + 2 \mbox{ match} (U - V')$.
\end{ex}

Note that if $U$ is an $n$-sun, $\smr (U) = 2 n$ also follows from Theorem~\ref{upf}, since in this case $U$ has a unique perfect matching. 

\begin{defn}
A pendant star, $S_\sharp$, in a unicyclic graph $U$, is an induced subgraph of $U$ such that
\begin{enumerate} 
\item $S_\sharp$ has exactly one vertex $v \in V_U$ with $\deg (v) \ge 2$,
\item $U - v$ has $k+1$ components, exactly $k$ of them are pendant vertices and the other is a (connected) unicyclic induced subgraph of $U$, and
\item $S_\sharp$ is induced by $v$ and the $k$ pendant vertices.
 \end{enumerate}
\end{defn}

\begin{alg}\label{algo}
Computation of $\smr (U)$, where $U$ is a connected unicyclic graph.
\end{alg}

Initialize: $U_\sharp = U$, $s = 0$

\begin{enumerate}

\item WHILE $U_\sharp$ is connected unicyclic and has a pendant star $S_\sharp$

\begin{enumerate}
\item $s=s+2$
\item $U_\sharp = U_\sharp- S_\sharp$
\end{enumerate}

END WHILE

\item If $U_\sharp$ is a cycle, compute $\smr (U_\sharp)$ as in Proposition~\ref{smr(C)}. If $U_\sharp$ is a partial $n$-dandelion, compute $\smr (U_\sharp)$ as in Example~\ref{pnsun}

\item $\smr (U) = \smr (U_\sharp) + s$
\end{enumerate}

\bpf
If $U_\sharp$ has a pendant star $S_\sharp$, then $v$  the center of $S_\sharp$, is a cut-vertex of $U_\sharp$, and $\smr (U_\sharp) = \smr (U_\sharp -v) + 2$. So we can let $U_\sharp = U_\sharp- S_\sharp$, and continue the process, inductively. Eventually, we are left with $U_\sharp$ as the unique cycle of the graph, or $U_\sharp$ as a partial $n$-dandelion, as in Definition~\ref{pndan}.
\epf

 We refer to the graph $U_\sharp$ produced by Algorithm~\ref{algo} as the {\em $\star$-reduced form} of $U$, and denote it by $U_\star$. The following theorem follows from Algorithm~\ref{algo}.

\begin{cor}\label{smruni}
If $U$ is connected unicyclic graph, then $$\smr (U) = 
\left\{
\begin{array}{ll}
2 \match (U) - 2 & \mbox{if $U_\star = C_k$ with $k$ even, $k \ge 4$},\\
2 \match (U) & \mbox{otherwise}.\\
\end{array}
\right.$$ 
\end{cor}

The following follows from Theorem~\ref{MR(G)}.
\begin{cor}
If $U$ is a unicyclic graph with unique cycle $C_k$, then $\smr (U) = \sMR (U)$ if and only if $k$ is odd, or $U_\star \ne C_k$ if $k$ even, $k \ge 4$.
\end{cor}

\begin{obs}\label{floor}
If $U$ is a connected unicyclic graph, with unique cycle $C = C_k$ and diametrical path $P$, then $|V_C \cap V_P| \leq \lfloor \frac{k}{2} \rfloor + 1$.
\end{obs}

In Lemma~\ref{necessary} we give necessary conditions on $U$ for $\smr (U) = \smr (P)$, where $P$ is a diametrical path of $U$. For $|P|$ odd we give necessary and sufficient conditions  in Theorem~\ref{mainuniodd}, the case $|P|$ even appears as Theorem~\ref{mainunieven}.
 
\begin{lem}\label{necessary}
Let $U$ be a connected unicyclic graph with unique cycle $C$, and diametrical path $P$. If $\smr (U) = \smr (P)$, then

\begin{enumerate}
\item\label{deg2} if  $v \in  V_C - V_P$, then $\deg_U (v) =2$;
\item\label{pede} every vertex $w \in V_C - V_P$ has the property that $U - w$ is a centipede;
\item\label{okcycles} the unique cycle is $C = C_k$, with $k \in \{3, 4, 6 \}$,
\begin{enumerate}
\item\label{onoff3}  if $k=3$, then $|V_{C_3} \cap V_P| = 2$ and $|V_{C_3} - V_P| = 1$;
\item\label{onoff4}  if $k=4$, then either $|V_{C_4} \cap V_P| = 2$ and $|V_{C_4} - V_P| = 2$, or  $|V_{C_4} \cap V_P| = 3$ and $|V_{C_4} - V_P| = 1$;
\item\label{onoff6}  if $k=6$, then $|V_{C_6} \cap V_P| = 4$ and $|V_{C_6} - V_P| = 2$.
\end{enumerate}
\end{enumerate}
\end{lem}

\bpf
\begin{enumerate}

\item If $v \in  V_C - V_P$ and $\deg_U (v) > 2$, then $U$ has an induced subgraph $U'$ with $ \smr (U') > \smr (P)$.

\item If $H = \{ w \} \subseteq V_C - V_P$, and $P'$ is  a diametrical path of $U-H$, then by Lemma~\ref{fallat}, $\smr (U - H) = \smr (U) = \smr (P) = \smr (P')$. By~\ref{deg2} above, $U-H$ is a tree, thus by Lemma~\ref{centi}, $U-H$ is a centipede.

\item  It is easy to see that  $|V_C \cap V_P| \ge 2$ and $|V_{C_k} - V_P| \le 2$. If $C_k$ is the unique cycle of $U$, then by Observation~\ref{floor}, we have $2 \le |V_{C_k} \cap V_P| \le \lfloor \frac{k}{2} \rfloor + 1$, which implies $3 \le k \le 6$, and $|V_C \cap V_P| \le 4$. Note that $k = 5$ is not a possibility. The three properties~\ref{onoff3},~\ref{onoff4} and~\ref{onoff6} now follow.

\end{enumerate}
\epf

\begin{thm}\label{mainuniodd}
Let $U$ be a connected unicyclic graph with unique cycle $C$ and diametrical path $P_n$, where $n \ge 5$ is odd. Then $\smr (U) = \smr (P_n)$ if and only if  $C = C_4$, with $|V_{C_4} \cap V_{P_n}| = 3$ and $|V_{C_4} - V_{P_n}| = 1$, the unique vertex $w \in V_{C_4} - V_{P_n}$ has the property that $U-w$ is a regular centipede with diametrical path $P_n$, and all joints are even-indexed vertices.
\end{thm}

\bpf
Let $P_n$, $n \ge 5$ odd, be a diametrical path of a regular centipede $T$, such that all joints are even-indexed vertices. Choose consecutive vertices $v_j, v_{j+1}, v_{j+2} \in V_{P_n}$, $j$ even, $2 \le j \le n-3$. Add a vertex $w$ to $T$, adjacent to $v_j$ and $v_{j+2}$. The resulting graph is a connected unicyclic graph $U$ with the properties described above, and by Corollary~\ref{smruni}, $\smr (P_n) = \smr (U)$. 

Conversely, suppose that  $U$ is a connected unicyclic graph with unique cycle $C$, diametrical path $P_n$, where $n \ge 5$ is odd, and $\smr (U) = \smr (P_n)$.  Let the unique cycle in $U$ be $C_k$, and $U'$ be the subgraph of $U$ induced by $V_{C_k} \cup V_{P_n}$. By Corollary~\ref{smruni}, if $k  = 3$, $k = 4$ with $|V_{C_4} \cap V_{P_n}| = 2$ and $|V_{C_4} - V_{P_n}| = 2$, or $k=6$, $\smr(U)  >  \smr(P_n)$. Hence, the unique cycle is $C_4$, with $|V_{C_4} \cap V_{P_n}| = 3$ and $|V_{C_4} - V_{P_n}| = 1$. By~Lemma~\ref{necessary}, $U- w$ is a centipede, where $w$ is the unique vertex in $V_C - V_{P_n}$, and thus $P_n$ is a diametrical path of $U -w$. By~Theorem~\ref{maintree}, $U- w$ is a regular centipede and all joints are even-indexed vertices.
\epf

The next lemma provides useful examples, the proofs follow from Corollary~\ref{smruni}.

\begin{lem}\label{pineapple}
If  $P_n$ is a path with $n$ even, and $U_i, i = 1, 2, 3, 4$, as described below, then $\smr (U_i) > \smr (P_n)$.

\begin{enumerate}
\item\label{3oddl} Let $U_1$ be the connected unicyclic graph with vertex set $V_{P_n} \cup \{ w, z \}$ and edge set $E_{P_n} \cup \{ \{ v_j, w \}, \{ w, v_{j+1} \}, \{ v_i, z \} \}$, where $j$ is odd, $3 \le j \le n-1$, and $i$ is odd, $3 \le i \le j$. 
\item\label{3oddr} Let $U_2$ be the connected unicyclic graph with vertex set $V_{P_n} \cup \{ w, z \}$ and edge set $E_{P_n} \cup \{ \{ v_j, w \}, \{ w, v_{j+1} \}, \{ v_i, z \} \}$, where $j$ is odd, $1 \le j \le n-3$, and $i$ is even, $j+1 \le i \le n-2$. 
\item\label{3evenl} Let $U_3$ be the connected unicyclic graph with vertex set $V_{P_n} \cup \{ w, z \}$ and edge set $E_{P_n} \cup \{ \{ v_j, w \}, \{ w, v_{j+1} \}, \{ v_i, z \} \}$, where $j$ is even, $4 \le j \le n-2$, and $i$ is odd, $3 \le i \le j-1$. 
\item\label{3evenr} Let $U_4$ be the connected unicyclic graph with vertex set $V_{P_n} \cup \{ w, z \}$ and edge set $E_{P_n} \cup \{ \{ v_j, w \}, \{ w, v_{j+1} \}, \{ v_i, z \} \}$, where $j$ is even, $2 \le j \le n-4$, and $i$ is even, $j+2 \le i \le n-2$. 
\end{enumerate}

\end{lem}

\begin{figure}[h!]
\begin{center}
\scalebox{.5}{\includegraphics{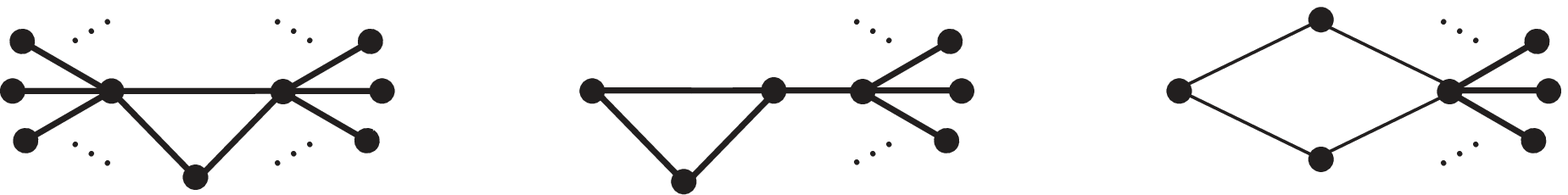}}
\caption{Small graphs for which $\smr (G) = \smr (P_4)$.}\label{small}
\end{center}
\end{figure}

For a connected unicyclic graph $U$, with unique cycle $C_k$, and diametrical path $P_n$, $n$ even, $n \le 4$, it is not difficult to see that $\smr (U) = \smr (P_n)$ if and only if $n=2$ and $U = C_3$, or $n=4$ and $U$ is one of the three graphs shown in Figure~\ref{small}, or  $U =C_6$.

In Lemma~\ref{conditions}, we give necessary conditions on a connected unicyclic graph $U$, for $\smr (U) = \smr (P_n)$, where $P_n$ is a diametrical path, with $n \ge 6$ and $n$ is even. Note that if  $T$ is a centipede with diametrical path $P_n=v_1 v_2 \cdots v_{n-1} v_n$, where $n$ is even, then reversing the numbering of the vertices changes odd-numbered vertices to even-numbered vertices and vice-versa; this action leads to an isomorphic centipede $T'$ with analogous properties. For the sake of completeness and to avoid confusion, we have included results assuming the original order of the vertices in $V_{P_n}$, although the symmetry of the graphs thus obtained introduces a slight redundancy. 

\begin{lem}\label{conditions}

Let $U$ be a connected unicyclic graph with unique cycle $C_k$ and diametrical path $P_n$, with $n$ even, $n \ge 6$, such that $\smr (U) = \smr (P_n)$.
\begin{enumerate}
\item\label{c3conditions}
If $k = 3$, then $V_{C_3} \cap V_{P_n} = \{ v_j, v_{j+1} \}$, and 
\begin{enumerate}
\item $1 \le j \le n-1$;
\item if $j=1$, then $\deg_U (v_1) = 2$, and $\deg_U (v_2) = 3$; 
\item if $j=n-1$, then $\deg_U (v_n) = 2$, and $\deg_U (v_{n-1}) = 3$; 
\item if $j$ is odd, $3 \le j \le n-3$, then  $\deg_U (v_j) = \deg_U(v_{j+1}) =3$;
\item for $i$ odd, $1 < i < j$, $\deg_U (v_i) = 2$, and for $i$ even, $ j+1 < i < n$, $\deg_U (v_i) = 2$.
\end{enumerate}
\item\label{c4conditions1} If $k=4$, and $V_{C_4} \cap V_{P_n} = \{ v_j, v_{j+1}, v_{j+2} \}$, then
\begin{enumerate}
\item\label{4a} $1 \le j \le n-2$;
\item\label{4b} $\deg_U (v_{j+1}) = 2$;
\item\label{4c} if $j = 1$, then $\deg_U (v_1) = 2$;
\item\label{4d}  if $j = n-2$, then $\deg_U (v_n) = 2$;
\item\label{4e} if $j=1$ and $\deg_U (v_3) > 3$, then for $i$ even, $3 < i < n$, $\deg_U (v_i) = 2$;
\item\label{4f} if $j=n$ and $\deg_U (v_{n-2}) > 3$, then for $i$ odd, $1 < i < n-2$, $\deg_U (v_i) = 2$;
\item\label{4g} if $j$ is odd, $3 \le j \le n-3$, and $\deg_U (v_j) > 3$ or $\deg_U (v_{j+2}) > 3$, then for $i$ even, $j + 2 < i < n$, $\deg_U (v_i) = 2$;
\item\label{4h} if $j$ is even, $2 \le j \le n-4$, and $\deg_U (v_j) > 3$ or $\deg_U (v_{j+2}) > 3$, then for $i$ odd, $1 < i < j$, $\deg_U (v_i) = 2$;
\item\label{4i} if for some $r$ odd, $r \ne j$, $r \ne j+2$, and $3 \le r \le n-3$, $\deg_U (v_r) > 2$, then for $i$ even, $r< i < n, i \ne j, i \ne j+2$, $\deg_U (v_i) = 2$;
\item\label{4j} if for some $r$ even, $r \ne j$, $r \ne j+2$, and $4 \le r \le n-2$, $\deg_U (v_r) > 2$, then for $i$ odd, $1 < i < r, i \ne j, i \ne j+2$, $\deg_U (v_i) = 2$.
\end{enumerate}
\item\label{c4conditions2} If $k=4$, and $V_{C_4} \cap V_{P_n} = \{ v_j, v_{j+1} \}$, then
\begin{enumerate}
\item  $j$ is odd, with $3 \le j \le n-3$;
\item  $\deg_U (v_j) = \deg_U (v_{j+1}) = 3$;
\item for $i$ odd, $1 < i < j$, $\deg_U (v_i) = 2$, and for $i$ even, $ j+1 < i < n$, $\deg_U (v_i) = 2$.
\end{enumerate}
\item\label{c6conditions} If $k=6$, then  $V_{C_6} \cap V_{P_n} = \{ v_j, v_{j+1}, v_{j+2}, v_{j+3} \}$, and 
\begin{enumerate}
\item $j$ is odd, with $1 \le j \le n-3$; 
\item $\deg_U (v_{j+1}) =  \deg_U (v_{j+2}) =2$;
\item if $j=1$, then $\deg_U (v_1) = 2$, and $\deg_U (v_4) = 3$;
\item\label{jn3} if $j=n-3$, then $\deg_U (v_n) = 2$, and $\deg_U (v_{n-3}) = 3$; 
\item if $j$ is odd, $3 \le j \le n-5$, then  $\deg_U (v_j) = \deg_U(v_{j+3}) =3$;
\item for $i$ odd, $1 < i < j$, $\deg_U (v_i) = 2$, and for $i$ even, $ j+1 < i < n$, $\deg_U (v_i) = 2$.
\end{enumerate} 
\end{enumerate}
\end{lem}

\bpf
\begin{enumerate}
\item Since $P_n$ is a diametrical path, if $j=1$, then $\deg_U (v_1) = 2$. If $\deg_U (v_2) > 3$, then taking  $i=2$, $U$ has $U_{\ref{3oddr}}$ as an induced subgraph (The graphs $U_i, i = 1,2,3,4$ refer to the graphs in Lemma~\ref{pineapple}). The case $j=n-1$ follows by letting $i = n-1$, and using the graph $U_{\ref{3oddl}}$. If $j$ is odd, with $3 \le j \le n-3$, and either $\deg_U (v_j) > 3$, or $\deg_U (v_{j+1}) > 3$, then letting $i=j$, or $i = j+1$, respectively, $U$ has as an induced subgraph, either the graph $U_{\ref{3oddl}}$, or the graph $U_{\ref{3oddr}}$, respectively. If $i$ is odd, $1 < i < j$, and $\deg_U (v_i) > 2$, then $U$ has either the graph $U_{\ref{3oddl}}$, or the graph $U_{\ref{3evenl}}$ as an induced subgraph. The case $j+1 < i < n$, with $i$ even follows by using either the graph $U_{\ref{3oddr}}$, or the graph $U_{\ref{3evenr}}$. 

\item Clearly $\deg_U (v_1) = 2$,  if $j=1$, and $\deg_U (v_n) = 2$ if $j = n-2$. If $w$ is the unique vertex in $V_C - V_{P_n}$, then the path $P' = v_1 v_2 \cdots v_j w v_{j+2} \cdots v_{n-1} v_n$, is also a diametrical path of $U$, with $v_{j+1} \in V_{C_4} - V_{P'}$, and thus  from Lemmas~\ref{fallat} and ~\ref{necessary}, $\deg_U (v_{j+1}) = 2$, $\smr (U) = \smr (U-w)$, $U-w$ is a centipede and has diametrical path $P_n$, and hence by Theorem~\ref{maintree}, properties~\ref{4e},~\ref{4f},~\ref{4g},~\ref{4h},~\ref{4i}, and~\ref{4j} above, must hold.

\item If $V_{C_4} \cap V_{P_n} = \{ v_j, v_{j+1} \}$, with $j$ even, $2 \le j \le n-2$, then by Corollary~\ref{smruni}, $\smr (U) > \smr (P_n)$. Thus $j$ must be odd, and clearly $3 \le j \le n-3$. Let $V_{C_4} - V_{P_n} = \{ w_1, w_2 \}$, with $w_1$ adjacent to $v_j$ and $w_2$ adjacent to $v_{j+1}$. By  Lemmas~\ref{fallat} and~\ref{necessary}, $\smr (U) = \smr (U - w_1) = \smr (U - w_2) = \smr (P_n)$, and both $U - w_1$ and $U - w_2$ are centipedes. Also, since $P_n$ is a diametrical path of either $U - w_1$ or $U - w_2$, by Theorem~\ref{maintree} no pair of  joints $v_r$ and $v_s$ satisfy $r$ odd, $s$ even, and $3 \le r < s \le n-2$. Since in the centipede $U- w_1$, $v_{j+1}$ is a joint, and $j+1$ is even, then $\deg_U (v_j) = 3$, and for $i$ odd, $1 < i < j$, $\deg_U (v_i) = 2$. Similarly $\deg_U (v_{j+1}) = 3$, and for $i$ even, $j+1 < i < n$, $\deg_U (v_i) = 2$.

\item If   $V_{C_6} \cap V_{P_n} = \{ v_j, v_{j+1}, v_{j+2}, v_{j+3} \}$, with $2 \le j \le n-3$, $j$ even,  by Corollary~\ref{smruni} $\smr (U) > \smr (P_n)$. Let $w_1, w_2 \in V_C - V_{P_n}$ with $w_1$ adjacent to $v_j$, $w_2$ adjacent to $v_{j+3}$, and $w_1$ adjacent to $w_2$. Then $U$ has a diametrical path the  path $P'' = v_1 v_2 \cdots v_j w_1 w_2 v_{j+3} \cdots v_{n-1} v_n$, thus by Lemma~\ref{necessary}, $\deg_U (v_{j+1}) = \deg_U (v_{j+2}) = 2$. If $j=1$, then clearly, $\deg_U (v_1) = 2$. Also, a diametrical path of $U-w_2$ is now the path $P' = w_1 v_1 v_2 \cdots v_{n-1} v_n$, and $\smr (U) = \smr (P')$. By Lemma~\ref{necessary}, $U-w_2$ is a centipede, thus by Theorem~\ref{maintree}, all the joints in $P'$ are odd-indexed vertices of $P_n$, therefore $\deg_U (v_4) = 3$. The proof of~\ref{jn3} is similar. Now let $j$ be odd, $3 \le j \le n - 5$. By Lemmas~\ref{fallat} and~\ref{necessary}, $\smr (U) = \smr (U-w_1) = \smr (P_n)$, $U-w_1$ is a centipede and has diametrical path $P_n$, thus by Theorem~\ref{maintree}, no pair of  joints $v_r$ and $v_s$ satisfy $r$ odd, $s$ even, with $3 \le r < s \le n-2$. But $U-w_1$ has a joint at $v_{j+3}$, which is an even-indexed vertex, that means $v_j$ cannot be a joint, and neither can any vertex $v_i$ with $i$ odd, $1 < i < j$. This implies $\deg_U (v_j) = 3$ and $\deg_U (v_i) = 2$, for $i$ odd, $1 < i < j$. Similarly, looking at $U-w_2$, we arrive at $\deg_U (v_{j+3}) = 3$ and for $i$ even, $ j+1 < i < n$, $\deg_U (v_i) = 2$.
\end{enumerate}
\epf

\begin{thm}\label{mainunieven}
Let $U$ be a connected unicyclic graph with unique cycle $C_k$ and diametrical path $P_n$, with $n$ even $n \ge 6$. Then $\smr (U) = \smr (P_n)$ if and only if
\begin{enumerate}
\item\label{path} every  vertex $w \in V_{C_k} - V_{P_n}$ has the property that $U-w$ is a centipede with diametrical path $P_n$, such that no pair of  joints $v_r$ and $v_s$ satisfy $r$ odd, $s$ even, and $3 \le r < s \le n-2$, or $U-w$ is a centipede with diametrical path $P'$ with vertex set $V_{P_n} \cup \{ z \}$, where $z$ is another vertex in $V_{C_k} - V_{P_n}$, and is adjacent to either $v_1$ or $v_n$,  such that all  joints of  $P'$ are odd-indexed vertices of $P_n$ whenever $z$ is adjacent to $v_1$, or all  joints of  $P'$ are even-indexed vertices in $P_n$ whenever $z$ is adjacent to $v_n$.

\item\label{cycles} $k = 3, 4$ or $6$.
\item\label{c3}
If $k = 3$, then $V_{C_3} \cap V_{P_n} = \{ v_j, v_{j+1} \}$, and 
\begin{enumerate}
\item $1 \le j \le n-1$;
\item if $j=1$, then $\deg_U (v_1) = 2$, and $\deg_U (v_2) = 3$; 
\item if $j=n-1$, then $\deg_U (v_n) = 2$, and $\deg_U (v_{n-1}) = 3$; 
\item if $j$ is odd, $3 \le j \le n-3$, then  $\deg_U (v_j) = \deg_U(v_{j+1}) =3$;
\item for $i$ odd, $1 < i < j$, $\deg_U (v_i) = 2$, and for $i$ even, $ j+1 < i < n$, $\deg_U (v_i) = 2$.
\end{enumerate}
\item\label{c4a} If $k=4$, and $V_{C_4} \cap V_{P_n} = \{ v_j, v_{j+1}, v_{j+2} \}$, then
\begin{enumerate}
\item\label{4a2} $1 \le j \le n-2$;
\item\label{4b2} $\deg_U (v_{j+1}) = 2$;
\item\label{4c2} if $j = 1$, then $\deg_U (v_1) = 2$;
\item\label{4d2}  if $j = n-2$, then $\deg_U (v_n) = 2$;
\item\label{4e2} if $j=1$ and $\deg_U (v_3) > 3$, then for $i$ even, $3 < i < n$, $\deg_U (v_i) = 2$;
\item\label{4f2} if $j=n$ and $\deg_U (v_{n-2}) > 3$, then for $i$ odd, $1 < i<n-2$, $\deg_U (v_i) = 2$;
\item\label{4g2} if $j$ is odd, $3 \le j \le n-3$, and $\deg_U (v_j) > 3$ or $\deg_U (v_{j+2}) > 3$, then for $i$ even, $j+2 < i < n$, $\deg_U (v_i) = 2$;
\item\label{4h2} if $j$ is even, $2 \le j \le n-4$, and $\deg_U (v_j) > 3$ or $\deg_U (v_{j+2}) > 3$, then for $i$ odd, $1 < i < j$, $\deg_U (v_i) = 2$;
\item\label{4i2} if for some $r$ odd, $r \ne j$, $r \ne j+2$, and $3 \le r \le n-3$, $\deg_U (v_r) > 2$, then for $i$ even, $r< i < n, i \ne j, i \ne j+2$, $\deg_U (v_i) = 2$;
\item\label{4j2} if for some $r$ even, $r \ne j$, $r \ne j+2$, and $4 \le r \le n-2$, $\deg_U (v_r) > 2$, then for $i$ odd, $1 < i < r, i \ne j, i \ne j+2$, $\deg_U (v_i) = 2$.\end{enumerate}
\item\label{c4b} If $k=4$, and $V_{C_4} \cap V_{P_n} = \{ v_j, v_{j+1} \}$, then
\begin{enumerate}
\item  $j$ is odd, with $3 \le j \le n-3$;
\item  $\deg_U (v_j) = \deg_U (v_{j+1}) = 3$;
\item for $i$ odd, $1 < i < j$, $\deg_U (v_i) = 2$, and for $i$ even, $ j+1 < i < n$, $\deg_U (v_i) = 2$.
\end{enumerate}
\item\label{c6} If $k=6$, then  $V_{C_6} \cap V_{P_n} = \{ v_j, v_{j+1}, v_{j+2}, v_{j+3} \}$, and 
\begin{enumerate}
\item $j$ is odd, with $1 \le j \le n-3$; 
\item $\deg_U (v_{j+1}) =  \deg_U (v_{j+2}) =2$;
\item if $j=1$, then $\deg_U (v_1) = 2$, and $\deg_U (v_4) = 3$;
\item if $j=n-3$, then $\deg_U (v_n) = 2$, and $\deg_U (v_{n-3}) = 3$; 
\item if $j$ is odd, $3 \le j \le n-5$, then  $\deg_U (v_j) = \deg_U(v_{j+3}) =3$;
\item for $i$ odd, $1 < i < j$, $\deg_U (v_i) = 2$, and for $i$ even, $ j+1 < i < n$, $\deg_U (v_i) = 2$.
\end{enumerate} 
\end{enumerate}
\end{thm}

\bpf
Let $U$ be a connected unicyclic graph with unique cycle $C$ and diametrical path $P_n$, $n$ even $n \ge 6$. Let $w \in V_C - V_{P_n}$, and assume that $U-w$ is a centipede.

Suppose $P_n$ is a diametrical path of $U-w$ with the property that no pair of  joints $v_r$ and $v_s$ satisfy $r$ odd, $s$ even, and $3 \le r < s \le n-2$. From Theorem~\ref{maintree} and Theorem~\ref{smrtree} it follows that $\smr (U-w) = \smr (P_n) =  2 \match (P_n)$.

If $k=3$, and $U$ satisfies the properties in~\ref{c3}, or $k=4$, $V_{C_4} \cap V_{P_n} = \{ v_j, v_{j+1}, v_{j+2} \}$, and $U$ satisfies the properties in~\ref{c4a}, then from Corollary~\ref{smruni} it follows that $\smr (U) = \smr (P_n)$.

If $k=4$, $V_{C_4} \cap V_{P_n} = \{ v_j, v_{j+1} \}$ and $U$ satisfies the properties in~\ref{c4b}, or $k=6$, $V_{C_6} \cap V_{P_n} = \{ v_j, v_{j+1}, v_{j+2}, v_{j+3} \}$ and $U$ satisfies the properties in~\ref{c6}, and either  $j \ne 1$, and $j \ne n-3$, or $j=1$ and $w$ is adjacent to $v_1$, or $j=n-3$, and $w$ is adjacent to $v_n$, then from Corollary~\ref{smruni} it follows that $\smr (U) =  2 \match (U) - 2 = \smr (P_n)$.

Now suppose $U-w$ has diametrical path $P'$ with vertex set $V_{P_n} \cup \{ z \}$, where $z$ is another vertex in $V_C - V_{P_n}$ adjacent to $v_1$, and all  joints of  $P'$ are odd-indexed vertices of $P_n$.

If $k=6$, and $j=1$, with $w$ adjacent to $v_4$, and $U$ satisfies all the properties in~\ref{c6}, then from Corollary~\ref{smruni}, $\smr (U) =  2 \match (U) - 2 = \smr (P_n)$.

Similarly, suppose $U-w$ has diametrical path $P'$ with vertex set $V_{P_n} \cup \{ z \}$, where $z$ is another vertex in $V_C - V_{P_n}$ adjacent to $v_n$, and all  joints of  $P'$ are even-indexed vertices in $V_{P_n}$. If $k=6$, and $j=n-3$, with $w$ adjacent to $v_{n-3}$, and $U$ satisfies all the properties in~\ref{c6}, then $\smr (U) = \smr (P_n)$.

Conversely, let $U$ be a connected unicyclic graph with unique cycle $C$ and diametrical path $P_n$, and $n$ even, $n \ge 6$, and let $\smr (U) = \smr (P_n)$. By Lemma~\ref{necessary}, the unique cycle is $C_k$, with $k \in \{3, 4, 6 \}$ (so condition~\ref{cycles} is met), and every vertex $w \in V_C - V_{P_n}$ has the property that $U - w$ is a centipede. By Lemmas~\ref{necessary}, and~\ref{conditions}, conditions~\ref{c3},~\ref{c4a},~\ref{c4b}, and~\ref{c6}, above are met. 

If $k = 3$, or $k = 4$, and $w \in V_C - V_{P_n}$, then $U-w$ is a centipede with diametrical path $P_n$. By Lemma~\ref{fallat}, $\smr (U - w) = \smr (U) = \smr (P_n)$, hence by~Theorem~\ref{maintree}, $P_n$ has the property that no pair of  joints $v_r$ and $v_s$ satisfy $r$ odd, $s$ even, with $3 \le r < s \le n-2$, thus condition~\ref{path} above is met. 

If $k = 6$, and $w \in V_C - V_{P_n}$, $w$ is not adjacent to $v_4$ or $v_{n-3}$, then $U-w$ is a centipede with diametrical path $P_n$. By Lemma~\ref{fallat}, $\smr (U - w) = \smr (U) = \smr (P_n)$, hence by~Theorem~\ref{maintree}, $P_n$ has the property that no pair of  joints $v_r$ and $v_s$ satisfy $r$ odd, $s$ even, with $3 \le r < s \le n-2$, thus condition~\ref{path} above is met. If $w$ is adjacent to either $v_4$ or $v_n-3$, then $U-w$ is a centipede with diametrical path $P'$ with vertex set $V_{P_n} \cup \{ z \}$, where $z$ is the other vertex in $V_C - V_{P_n}$. By Lemma~\ref{fallat}, $\smr (U - w) = \smr (U) = \smr (P_n)$, by~Theorem~\ref{maintree}, $P'$ has the property that all joints are odd-indexed vertices of $P_n$, whenever $z$ is adjacent to $v_1$, or all joints of $P'$ are even-indexed vertices in $P_n$, whenever $z$ is adjacent to $v_n$, thus condition~\ref{path} above is met.
\epf

\section{Additional Examples}\label{examples}

\begin{figure}[h!]
\begin{center}
\scalebox{.5}{\includegraphics{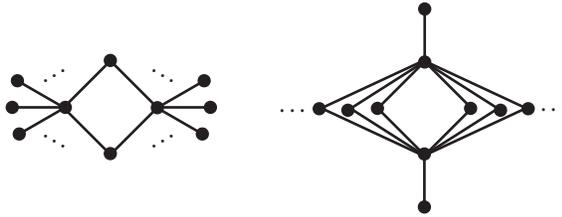}}
\caption{A spider $S_k$ and a butterfly $B_k$, graphs for which $\smr (G) = \dm (G)$.}\label{bugs}
\end{center}
\end{figure}

\begin{ex}\label{uni4}  The connected graphs $S_k$ and $B_k$, in Figure~\ref{bugs} satisfy $\smr (S_k) = \dm(S_k) = 4$, and $\smr (B_k) = \dm (B_k) = 4$.
\end{ex}

The spider $S_k$,  clearly satisfies the conditions of Theorem~\ref{mainuniodd}, and $\smr (B_k) = \dm (B_k)$ follows from cut-vertex reduction, 

Finally we show some general types of connected not unicyclic graphs $G$, with the property that $\smr (G) = \dm (G)$.

\begin{figure}[h!]
\begin{center}
\scalebox{.5}{\includegraphics{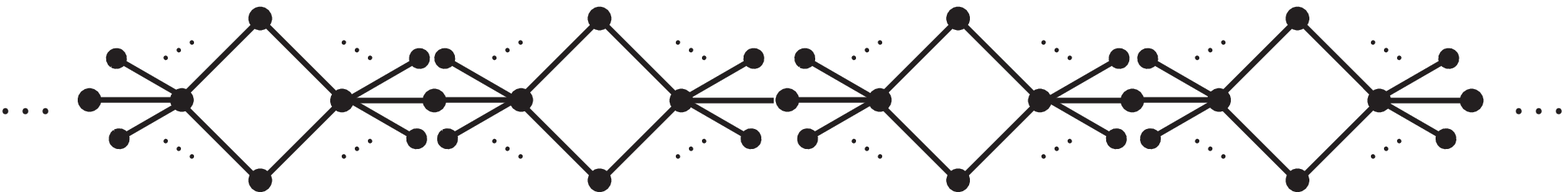}}
\caption{A connected graph for which $\smr (G) = \dm (G)$.}\label{spiders}
\end{center}
\end{figure}

\begin{prop}
The graph shown in Figure~\ref{spiders} satisfies $\smr (G) = \dm (G)$.
\end{prop}

\bpf
The proof is by induction on the number $n$, of spiders in $G$. Let $G_s$ denote a graph, as pictured in Figure~\ref{spiders}, consisting of exactly $s$  of these spiders; also, note that the path $P_{4 s +1}$ is an induced subgraph of $G_s$, and that $\dm (G_s) = 4s$. The base case, $n=1$, follows from Example~\ref{uni4}.  Assume the graph, $G_k$ satisfies  $\smr (G_k) = \dm (G_k)$, and let $G_{k+1}$ be the graph, as in Figure~\ref{spiders},  with one more spider than $G_k$. The graph $G_{k+1}$ is the vertex sum of $G_k$ and $G_1$ at the vertex $v$, of degree 2, between the spiders. Since $P_4$ is an induced path in $G_1 - v$, it follows that $\smr (G_1) = \smr (G_1 - v)$, and thus $\srs_v (G_1) = 0$. From Theorem~\ref{cutvertex}, we have $\smr (G_{k+1}) = \smr (G_k - v) + \srs_v (G_k) + \smr (G_1 - v) = \smr (G_k) + \smr (G_1) = \dm (G_k) + 4 = 4k + 4 =4 (k + 1)$, as we wanted to show. 
\epf

\begin{figure}[h!]
\begin{center}
\scalebox{.5}{\includegraphics{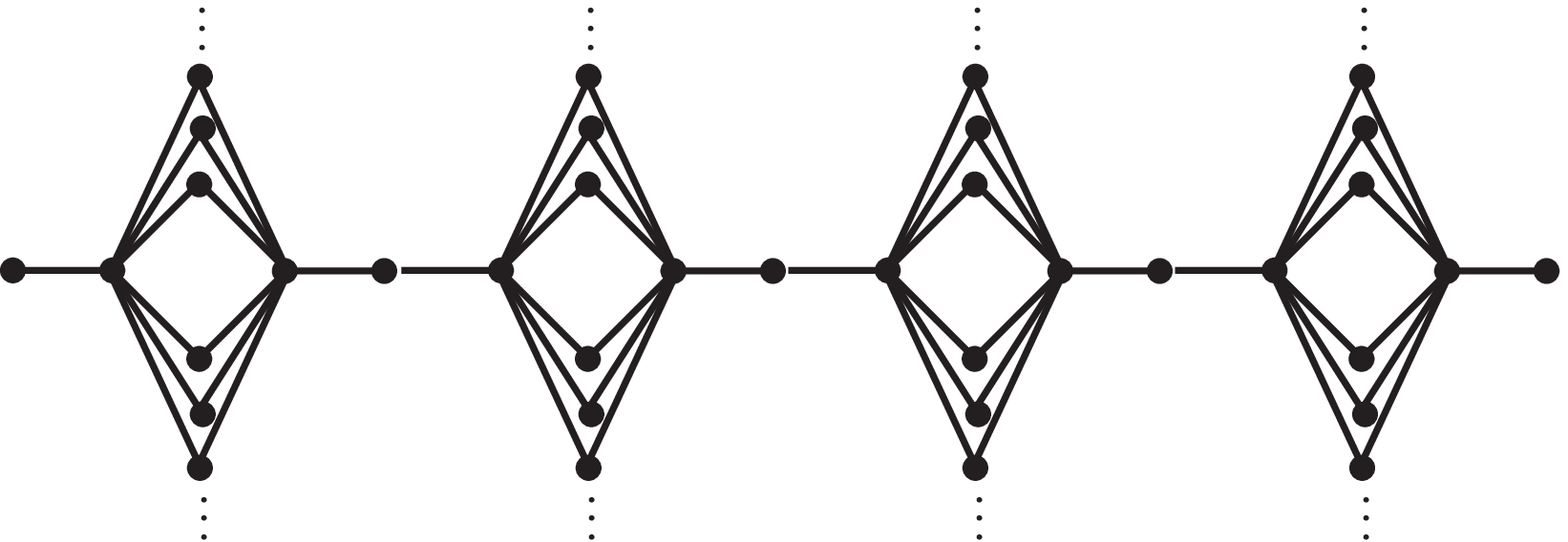}}
\caption{A connected graph for which $\smr (G) = \dm (G)$.}\label{bflies}
\end{center}
\end{figure}

\begin{figure}[h!]
\begin{center}
\scalebox{.5}{\includegraphics{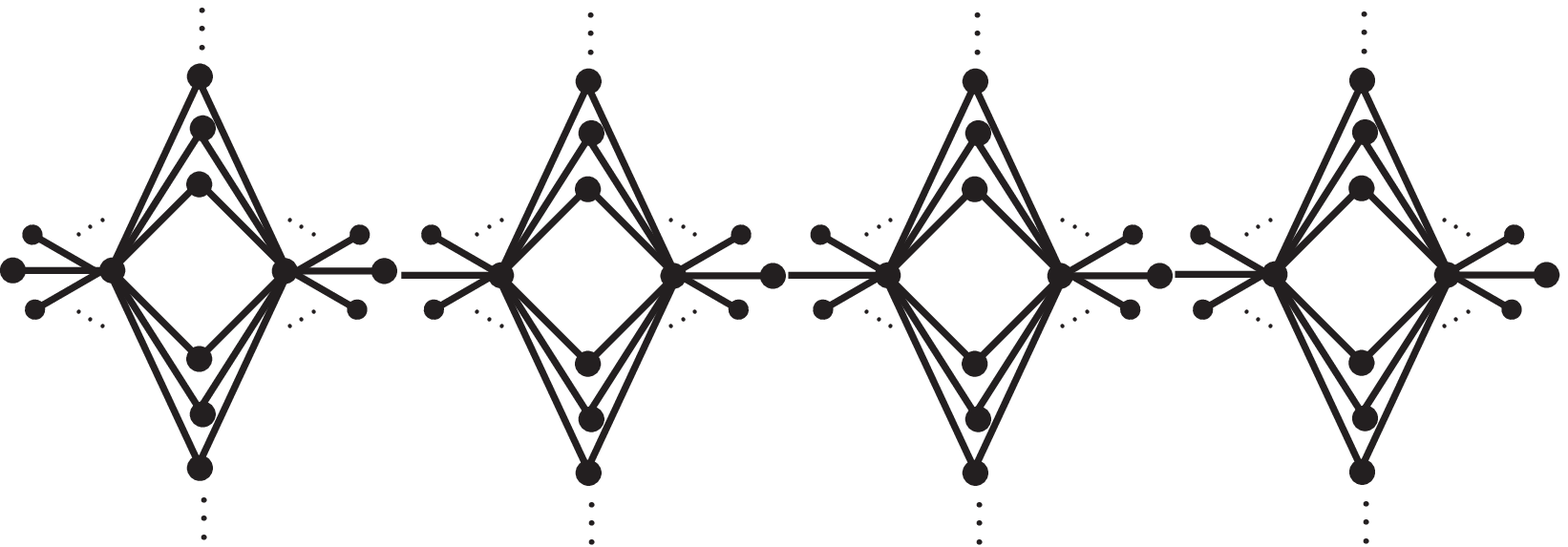}}
\caption{A connected graph for which $\smr (G) = \dm (G)$.}\label{allbugs}
\end{center}
\end{figure}

\begin{prop}
Induction arguments similar to the one in Proposition~\ref{spiders} show that the graphs shown in Figure~\ref{bflies}, and Figure~\ref{allbugs} satisfy $\smr (G) = \dm (G)$.
\end{prop}

{\bf Acknowledgement.} The author wishes to thank Bryan  L. Shader for suggesting the minimum skew rank problem at the ``2008 IMA PI Summer Program:
Linear Algebra and Applications.'' Iowa State University, Ames,
Iowa.  June  30--July 25, 2008.

\end{document}